# On q-series identities with multiple sums related to divisor functions

*Aung Phone Maw*

**Abstract**

We provide an exposition of q-identities with multiple sums related to divisor functions given by Dilcher, Prodinger, Fu and Lascoux, Zeng, Guo and Zhang. Meanwhile, for each of these identities, a more powerful statement will be derived through our exposition.

## 1. Introduction

Throughout this paper, we shall use the following standard notation

$$(x;q)_\infty = (x)_\infty = \prod_{i\geq 1}\left(1 - xq^{i-1}\right),$$

$$(x;q)_n = \frac{(x;q)_\infty}{(xq^n;q)_\infty}, n \in \mathbb{C}, \qquad (x;q)_0 = 1\,.$$

The study of *q-identities related to divisor functions*[1,2,4,9,11-14,16] has given rise to numerous interesting q-identities. In these studies, the regular appearance of multiple sums is of noticeable significance. We shall now mention those identities with multiple sums which will be examined throughout our study. The first appearance of these type of identities occurs in [4], where Dilcher gave the following

$$\sum_{n\geq i_1\geq\cdots\geq i_m\geq 1}\prod_{1\leq j\leq m}\frac{q^{i_j}}{1-q^{i_j}} = \sum_{1\leq r\leq n}\begin{bmatrix} n\\ r\end{bmatrix}\frac{(-1)^{r-1}q^{\binom{r}{2}+mr}}{(1-q^r)^m}, \tag{1.1}$$

where $\begin{bmatrix} n\\ r\end{bmatrix} = \frac{[n][n-1]\ldots[n-r+1]}{[r][r-1]\ldots[1]}$ is the *Gaussian binomial coefficient*. $[n] = \frac{1-q^n}{1-q}$ being the *q-number*. (1.1) is provided as a certain analogue of the identity

$$\sum_{i\geq 1}\frac{q^i}{1-q^i} = \sum_{r\geq 1}\frac{(-1)^{r-1}q^{\binom{r+1}{2}}}{(q;q)_r(1-q^r)}, \tag{1.2}$$

given in [5]. Later, Prodinger[1] proved the following

$$\sum_{1\leq i\leq n}\frac{q^{i(m-1)}}{(1-q^i)^m} = \sum_{1\leq r\leq n}\begin{bmatrix} n\\ r\end{bmatrix}(-1)^{r-1}q^{\binom{r}{2}-rn}\sum_{\substack{i_1\geq\cdots\geq i_m\geq 1\\ i_1=r}}\prod_{1\leq j\leq m}\frac{q^{i_j}}{1-q^{i_j}}, \tag{1.3}$$

by inverting the original result of Dilcher and thus giving a q-analog of a formula of Hernández[3]. And Fu and Lascoux[12] generalized (1.1) as

$$\sum_{n\geq i_1\geq\cdots\geq i_m\geq 1}(-1)^{i_m-1}\left(x^{i_m}-(-1)^{i_m}\right)\prod_{1\leq j\leq m}\frac{q^{i_j}}{1-q^{i_j}}=\sum_{1\leq r\leq n}\begin{bmatrix}n\\r\end{bmatrix}\frac{(-1)^{r-1}x^r(-x^{-1};q)_r q^{mr}}{(1-q^r)^m}\ . \tag{1.4}$$

For this, Fu and Lascoux used the Newton interpolation. While Prodinger[2] and Zeng[9] provided different proofs of (1.4). Zeng[9] in particular, using the method of partial fraction decomposition, obtained a further generalization of (1.4), which can be stated as

$$\frac{(q;q)_n}{(zq;q)_n}\sum_{n\geq i_1\geq\cdots\geq i_m\geq 1}\frac{x^{i_m}(zq;q)_{i_m}}{(q;q)_{i_m}}\prod_{1\leq j\leq m}\frac{q^{i_j}}{1-zq^{i_j}}=\sum_{1\leq r\leq n}\begin{bmatrix}n\\r\end{bmatrix}\frac{\left(x^r(x^{-1};q)_r+(-1)^{r-1}q^{\binom{r}{2}}\right)q^{mr}}{(1-zq^r)^m}\ . \tag{1.5}$$

Lately, this identity has become a very common generalization of (1.1). It is noteworthy that Ismail and Stanton[17], A.Xu[13] also provided different proofs of (1.5). And lastly, we mention Guo and Zhang[11], who obtained the following very unique generalization of (1.1)

$$-\sum_{i_1=1}^{n}\frac{q^{i_1}}{(1-q^{i_1})(1-zq^{i_1-1})}\sum_{i_2=1}^{i_1}\frac{q^{i_2}}{(1-q^{i_2})(1-zq^{i_2-2})}\cdots\sum_{i_m=1}^{i_{m-1}}\frac{q^{i_m}}{(1-q^{i_m})(1-zq^{i_m-m})}$$

$$=\sum_{1\leq r\leq n}\begin{bmatrix}n\\r\end{bmatrix}\frac{(z^{-1}q^m;q)_r(z;q)_{n-r}}{(zq^{-m};q)_{m+n}(1-q^r)^m}z^r\ . \tag{1.6}$$

In our study, we shall derive three main results. The first is a very general proposition generalizing $(1.1)-(1.4)$. The second is a proposition generalizing (1.5). And the last result generalizes (1.6). Before stating our main results, we recall that the $k^{th}$*complete symmetric function* of $a_1,a_2,\dots,a_n$, denoted as $h_k(a_1,a_2,\dots,a_n)$, is defined as

$$h_k(a_1,a_2,\dots,a_n)=\sum_{1\leq i_1\leq i_2\leq\cdots\leq i_k\leq n}a_{i_1}a_{i_2}\dots a_{i_k}\ ,$$

$$h_0(a_1,a_2,\dots,a_n)=1\ .$$

And the generating function of $h_k$ is

$$\sum_{i\geq 0}t^i h_i(a_1,a_2,\dots,a_n)=\frac{1}{(1-a_1t)(1-a_2t)\dots(1-a_nt)}\ . \tag{1.7}$$

Then, our main results are stated as follows.

**Theorem A.** *For* $m_j\geq 0, 1\leq j\leq k$, *there holds*

$$(a).\ \sum_{i_1=1}^{n}\frac{q^{m_1i_1}}{(1-q^{i_1})^{m_1+1}}\sum_{i_2=1}^{i_1}\frac{q^{m_2i_2}}{(1-q^{i_2})^{m_2+1}}\cdots\sum_{i_k=1}^{i_{k-1}}\frac{q^{m_ki_k}}{(1-q^{i_k})^{m_k+1}}x^{i_k}$$

$$=\sum_{i\geq 1}\begin{bmatrix}n\\i\end{bmatrix}\frac{(-1)^{i-1}q^{\binom{i+1}{2}-ni}}{1-q^i}\sum_{\substack{i\geq i_1\geq\cdots\geq i_k\geq 1\\ i_0=i}}\left(1-q^{i_k}\right)\left(1-(x;q)_{i_k}\right)\prod_{1\leq j\leq k}\frac{q^{i_j}}{\left(1-q^{i_j}\right)^2}h_{m_j-1}\left(\frac{q^{i_j}}{1-q^{i_j}},\dots,\frac{q^{i_{j-1}}}{1-q^{i_{j-1}}}\right), \tag{1.8}$$

$$(b).\quad \sum_{i_1=1}^{n}\frac{q^{i_1}}{(1-q^{i_1})^{m_1+1}}\sum_{i_2=1}^{i_1}\frac{q^{i_2}}{(1-q^{i_2})^{m_2+1}}\cdots\sum_{i_k=1}^{i_{k-1}}\frac{q^{i_k}}{(1-q^{i_k})^{m_k+1}}\left(1-(x;q)_{i_k}\right)$$

$$=\sum_{i\geq 1}\begin{bmatrix}n\\i\end{bmatrix}\frac{(-1)^{i-1}q^{\binom{i+1}{2}}}{1-q^i}\sum_{\substack{i\geq i_1\geq\cdots\geq i_k\geq 1\\ i_0=i}}\left(q^{-i_k}-1\right)x^{i_k}\prod_{1\leq j\leq k}\frac{q^{i_j}}{\left(1-q^{i_j}\right)^2}h_{m_j-1}\left(\frac{1}{1-q^{i_j}},\dots,\frac{1}{1-q^{i_{j-1}}}\right).\qquad(1.9)$$

In the context of this paper, it is plausible to treat the sum $\sum_{1\leq i\leq n}a_i\frac{q^i}{1-q^i}h_{m-1}\left(\frac{q^i}{1-q^i},\dots,\frac{q^n}{1-q^n}\right)$ as $a_n$ whenever $m=0$ . The same holds for the sum $\sum_{1\leq i\leq n}a_i\frac{1}{1-q^i}h_{m-1}\left(\frac{1}{1-q^i},\dots,\frac{1}{1-q^n}\right)$ (We shall see why in Section 2). Thus, for $m_1=\cdots=m_k=0,$ (1.8) results in (1.4), which is the result of Fu and Lascoux, since the sum

$$\sum_{\substack{i\geq i_1\geq\cdots\geq i_k\geq 1\\ i_0=i}}\left(1-q^{i_k}\right)\left(1-(x;q)_{i_k}\right)\prod_{1\leq j\leq k}\frac{q^{i_j}}{\left(1-q^{i_j}\right)^2}h_{m_j-1}\left(\frac{q^{i_j}}{1-q^{i_j}},\dots,\frac{q^{i_{j-1}}}{1-q^{i_{j-1}}}\right)$$

reduces to $\frac{1-(x;q)_i}{(1-q^i)^{k-1}}$. And it can be seen that for $k=1$ and $x=\ 1$ , (1.8) results in (1.3). Moreover, **Theorem A** can be compared with the duality identity(Theorem 1) appearing in Bradley[15] which does not contain the parameter $x$ .

**Proposition B.** *For k natural numbers $l_1,l_2,\dots,l_k$ and k arbitrary constants $z_1,z_2,\dots,z_k$, there holds*

$$(a).\quad \sum_{\substack{n\geq i_1\geq\cdots\geq i_k\geq 1\\ i_0=n}}\frac{(x;q)_{i_k}(z_kq;q)_{i_k}}{(q;q)_{i_k}}\prod_{1\leq j\leq k-1}\frac{\left(z_jq;q\right)_{i_j}}{\left(1-q^{i_j}\right)\left(z_{j+1}q;q\right)_{i_j}}\prod_{1\leq j\leq k}\frac{q^{i_j}}{\left(1-z_jq^{i_j}\right)}h_{l_j-1}\left(\frac{q^{i_j}}{1-z_jq^{i_j}},\dots,\frac{q^{i_{j-1}}}{1-z_jq^{i_{j-1}}}\right)$$

$$=\frac{(z_1q;q)_n}{(q;q)_n}\sum_{i\geq 1}\begin{bmatrix}n\\i\end{bmatrix}\frac{(-1)^{i-1}q^{\binom{i}{2}+il_1}}{(1-z_1q^i)^{l_1}}\sum_{\substack{i\geq i_1\geq\cdots\geq i_{k-1}\geq 1\\ i_0=i}}\left(1-x^{i_{k-1}}\right)\prod_{1\leq j\leq k-1}\frac{q^{i_jl_{j+1}}}{\left(1-q^{i_j}\right)\left(1-z_{j+1}q^{i_j}\right)^{l_{j+1}}},\qquad(1.10)$$

$(b).$

$$\sum_{\substack{n\geq i_1\geq\cdots\geq i_k\geq 1\\ i_0=n}}\frac{\left(1-x^{i_k}\right)z_k{}^{i_k}(z_k{}^{-1}q;q)_{i_k}}{(q;q)_{i_k}}\prod_{1\leq j\leq k-1}\frac{q^{i_j}z_j{}^{i_j}z_{j+1}{}^{-i_j}\left(z_j{}^{-1}q;q\right)_{i_j}}{\left(1-q^{i_j}\right)\left(z_{j+1}{}^{-1}q;q\right)_{i_j}}\prod_{1\leq j\leq k}\frac{1}{\left(z_j-q^{i_j}\right)}h_{l_j-1}\left(\frac{1}{z_j-q^{i_j}},\dots,\frac{1}{z_j-q^{i_{j-1}}}\right)$$

$$=\frac{z_1{}^n(z_1{}^{-1}q;q)_n}{(q;q)_n}\sum_{i\geq 1}\begin{bmatrix}n\\i\end{bmatrix}\frac{(-1)^{i-1}q^{\binom{i+1}{2}-ni}}{(z_1-q^i)^{l_1}}\sum_{\substack{i\geq i_1\geq\cdots\geq i_{k-1}\geq 1\\ i_0=i}}(x;q)_{i_{k-1}}\prod_{1\leq j\leq k-1}\frac{q^{i_j}}{\left(1-q^{i_j}\right)\left(z_{j+1}-q^{i_j}\right)^{l_{j+1}}}.\qquad(1.11)$$

It can be seen that **Proposition B** is a multiple generalization of the result of Zeng, since for $k=1$, (1.11) results in (1.5).

**Proposition C.** *For natural number k, there holds*

$$\sum_{i_1=1}^{n}\frac{q^{i_1}}{(1-q^{i_1})(1-zq^{i_1-1})}\sum_{i_2=1}^{i_1}\frac{q^{i_2}}{(1-q^{i_2})(1-zq^{i_2-2})}\cdots\sum_{i_k=1}^{i_{k-1}}\frac{q^{i_k}}{(1-q^{i_k})(1-zq^{i_k-k})}\sum_{r=1}^{i_k}\left(\frac{1}{1-zq^{r-k-1}}-\frac{1}{1-tq^{r-1}}\right)$$

$$=\sum_{r\geq 1}\begin{bmatrix}n\\r\end{bmatrix}\frac{(q;q)_{r-1}(tz^{-1}q^k;q)_r(z;q)_{n-r}}{(zq^{-k};q)_{n+k}(t;q)_r(1-q^r)^k}z^r\ . \tag{1.12}$$

If we multiply both sides of **Proposition C** with $1-t$ and set $t=1$, we arrive at (1.6), the result of Guo and Zhang. Moreover, it can be easily seen that Theorem 4.1 in [11] also by Guo and Zhang can be derived by making some elementary manipulations and substitution on **Proposition C**.

The rest of the paper is organized as follows. In Section 2, we shall prove **Theorem A** by repeated application of the Jackson integral. In Section 3, we shall prove **Proposition B** and draw some consequences. Finally, in Section 4, we shall prove **Proposition C**.

## 2. Proof of Theorem A

We denote by $D_q$ , the *q-derivative* of a function $f$

$$D_qf(x)=\big(D_qf\big)(x)=\frac{f(x)-f(xq)}{x-xq}\ . \tag{2.1}$$

The *Jackson definite integral* of the function $f$ is defined in [6] as

$$\int_0^x f(t)d_qt=(1-q)\sum_{n\geq 0}q^nxf(q^nx)\ . \tag{2.2}$$

The Jackson integral and the q-derivative are related by the *"fundamental theorem of quantum calculus"*[18,p.73], which implies

(a) If $D_qF=f$ and $F(x)$ is continuous at $x=0$, then

$$\int_0^x f(t)d_qt=F(x)-F(0)\,. \tag{2.3}$$

(b) For any function $f$

$$D_q\int_0^x f(t)d_qt=f(x)\,. \tag{2.4}$$

Then, it is easy to deduce

$$\int_0^x t^{n-1}d_qt=\frac{x^n}{[n]}\,, \tag{2.5}$$

$$\int_0^x (tyq;q)_{n-1}d_qt=-\frac{(xy;q)_n}{y[n]}\ . \tag{2.6}$$

Here, we define some operators using the Jackson integral for the purpose of the proof.

$$P_q f(x)=\int_0^x \frac{f(t)}{t} d_q t \ ,(P_q)^m f(x)=(P_q)^{m-1}\left(P_q f(x)\right),(P_q)^0 f=f \ ,$$

and

$$T_q f(x)=\int_0^x \frac{f(1)-f(t)}{1-t} d_q t \ ,(T_q)^m f(x)=(T_q)^{m-1}\left(T_q f(x)\right),(T_q)^0 f=f \ .$$

We shall also use the *q-shift operator* denoted as $\eta_x$ as in [8] defined by

$$\eta_x f(x)=f(xq)\ ,{\eta_x}^m f(x)={\eta_x}^{m-1}\big(\eta_x f(x)\big)\ ,{\eta_x}^0 f=f \ .$$

Then, our proof is based on the following lemma.

**Lemma 2.1.** *For* $\mathrm{P}_q$ *,* $\mathrm{T}_q$ *defined as above, the following holds*

$$(a).\ (P_q)^m x^n=\frac{x^n}{[n]^m}\ , \tag{2.7}$$

$$(b).\ (P_q)^m(1-(x;q)_n)=(1-q)^m \sum_{1\le i\le n}(1-(xq^{-m};q)_i)\frac{q^i}{1-q^i} h_{m-1}\left(\frac{q^i}{1-q^i},\dots,\frac{q^n}{1-q^n}\right), \tag{2.8}$$

$$(c).\ (T_q)^m x^n=(1-q)^m \sum_{1\le i\le n} x^i \frac{1}{1-q^i} h_{m-1}\left(\frac{1}{1-q^i},\dots,\frac{1}{1-q^n}\right), \tag{2.9}$$

$$(d).\ (T_q)^m(1-(x;q)_n)=\frac{1-(x;q)_n}{[n]^m}\ . \tag{2.10}$$

*Proof of Lemma 2.1.*

$$(a).\ \ (P_q)^m x^n=(P_q)^{m-1}\int_0^x t^{n-1} d_q t=(P_q)^{m-1}\frac{x^n}{[n]}=\ \frac{1}{[n]}\left(P_q\right)^{m-1} x^n=\dots=\frac{x^n}{[n]^m}\,.$$

$$(b).\ \ (P_q)^m(1-(x;q)_n)=(P_q)^{m-1}\int_0^x \frac{(1-(t;q)_n)}{t} d_q t$$

Since $$\frac{(1-(t;q)_n)}{t}=\sum_{1\le i\le n} q^{i-1}(t;q)_{i-1}\ , \tag{2.11}$$

$$(P_q)^m(1-(x;q)_n)=(P_q)^{m-1}\sum_{1\le i_1\le n} q^{i_1-1}\int_0^x (t;q)_{i_1-1} d_q t=\sum_{1\le i_1\le n}\frac{q^{i_1}}{[i_1]}(P_q)^{m-1}\left(1-\left({}^{x}\!/_{q}\,;q\right)_{i_1}\right)=\cdots$$

$$=\sum_{1\le i_1\le n}\frac{q^{i_1}}{[i_1]}\sum_{1\le i_2\le i_1}\frac{q^{i_2}}{[i_2]}\cdots\sum_{1\le i_m\le i_{m-1}}\frac{q^{i_m}}{[i_m]}\left(1-\left({}^{x}\!/_{q^m}\,;q\right)_{i_m}\right)$$

$$=(1-q)^m\sum_{1\le i\le n}(1-(xq^{-m};q)_i)\frac{q^i}{1-q^i} h_{m-1}\left(\frac{q^i}{1-q^i},\dots,\frac{q^n}{1-q^n}\right).$$

$(c)$. $(T_q)^m x^n = (T_q)^{m-1}\int_0^x \frac{1-t^n}{1-t} d_q t = \sum_{1\le i_1\le n} (T_q)^{m-1}\int_0^x t^{i_1-1} d_q t = \sum_{1\le i_1\le n} \frac{1}{[i_1]} (T_q)^{m-1} x^{i_1} = \cdots$

$$= \sum_{1\le i_1\le n} \frac{1}{[i_1]} \sum_{1\le i_2\le i_1} \frac{1}{[i_2]} \cdots \sum_{1\le i_m\le i_{m-1}} \frac{1}{[i_m]} x^{i_m}$$

$$= (1-q)^m \sum_{1\le i\le n} x^i \frac{1}{1-q^i} h_{m-1}\left(\frac{1}{1-q^i}, \dots, \frac{1}{1-q^n}\right).$$

$(d)$. $(T_q)^m (1-(x;q)_n) = (T_q)^{m-1}\int_0^x \frac{(t;q)_n}{1-t} d_q t = (T_q)^{m-1}\int_0^x (tq;q)_{n-1} d_q t = \frac{1}{[n]} (T_q)^{m-1} (1-(x;q)_n)$

$$= \cdots = \frac{1-(x;q)_n}{[n]^m}.$$

Note that in view of (2.8) and (2.9) , it is plausible to treat the sums $\sum_{1\le i\le n} a_i \frac{q^i}{1-q^i} h_{m-1}\left(\frac{q^i}{1-q^i}, \dots, \frac{q^n}{1-q^n}\right)$ and $\sum_{1\le i\le n} a_i \frac{1}{1-q^i} h_{m-1}\left(\frac{1}{1-q^i}, \dots, \frac{1}{1-q^n}\right)$ as $a_n$ whenever $m=0$. Now, we shall prove **Theorem A**.

*Proof of Theorem A.*

$(a)$.

Our starting point is the following inverted form of the q-binomial theorem.(Prodinger[1] also used this for their result)

$$x^n = \sum_{i\ge 1} \begin{bmatrix} n \\ i \end{bmatrix} (-1)^{i-1} q^{\binom{i+1}{2}-in} (1-(x;q)_i). \tag{2.12}$$

Then, we apply a series of operators $(\eta_x{}^{m_k} P_q{}^{m_k} T_q) \dots (\eta_x{}^{m_2} P_q{}^{m_2} T_q)(\eta_x{}^{m_1} P_q{}^{m_1} T_q)$ to both sides, we get

$L.H.S$
$$= (\eta_x{}^{m_k} P_q{}^{m_k} T_q) \dots (\eta_x{}^{m_2} P_q{}^{m_2} T_q)(\eta_x{}^{m_1} P_q{}^{m_1} T_q)\, x^n$$

$$= \sum_{1\le i_1\le n} \frac{q^{m_1 i_1}}{[i_1]^{m_1+1}} (\eta_x{}^{m_k} P_q{}^{m_k} T_q) \dots (\eta_x{}^{m_2} P_q{}^{m_2} T_q) x^{i_1} = \cdots$$

$$= \sum_{1\le i_1\le n} \frac{q^{m_1 i_1}}{[i_1]^{m_1+1}} \sum_{1\le i_2\le i_1} \frac{q^{m_2 i_2}}{[i_2]^{m_2+1}} \cdots \sum_{1\le i_k\le i_{k-1}} \frac{q^{m_k i_k}}{[i_k]^{m_k+1}} x^{i_k}$$

$$= (1-q)^{m_1+m_2+\cdots+m_2+k} \sum_{n\ge i_1\ge\cdots\ge i_k\ge 1} x^{i_k} \prod_{1\le j\le k} \frac{q^{m_j i_j}}{(1-q^{i_j})^{m_j+1}},$$

$R.H.S$

$$= \sum_{i\ge 1} \begin{bmatrix} n \\ i \end{bmatrix} (-1)^{i-1} q^{\binom{i+1}{2}-in} (\eta_x{}^{m_k} P_q{}^{m_k} T_q) \dots (\eta_x{}^{m_2} P_q{}^{m_2} T_q)(\eta_x{}^{m_1} P_q{}^{m_1} T_q)\, (1-(x;q)_i).$$

Since

$$\left(\eta_x{}^{m_k}P_q{}^{m_k}T_q\right)\ldots(\eta_x{}^{m_2}P_q{}^{m_2}T_q)(\eta_x{}^{m_1}P_q{}^{m_1}T_q)\,(1-(x;q)_i)$$

$$=\frac{(1-q)^{m_1+1}}{1-q^i}\sum_{1\le i_1\le i}\frac{q^{i_1}}{1-q^{i_1}}h_{m_1-1}\left(\frac{q^{i_1}}{1-q^{i_1}},\ldots,\frac{q^i}{1-q^i}\right)\left(\eta_x{}^{m_k}P_q{}^{m_k}T_q\right)\ldots\left(\eta_x{}^{m_2}P_q{}^{m_2}T_q\right)\left(1-(x;q)_{i_1}\right)$$

$$=\cdots=\frac{(1-q)^{m_1+m_2+\cdots+m_k+k}}{1-q^i}\sum_{\substack{i\ge i_1\ge\cdots\ge i_k\ge 1\\ i_0=i}}\left(1-q^{i_k}\right)\left(1-(x;q)_{i_k}\right)$$

$$\times\prod_{1\le j\le k}\frac{q^{i_j}}{\left(1-q^{i_j}\right)^2}h_{m_j-1}\left(\frac{q^{i_j}}{1-q^{i_j}},\ldots,\frac{q^{i_{j-1}}}{1-q^{i_{j-1}}}\right),$$

dividing both $L.H.S$ and $R.H.S$ by $(1-q)^{m_1+m_2+\cdots+m_k+k}$, we arrive at (1.8).

$(b)$.

For (1.9), our starting point is the q-binomial theorem itself

$$1-(x;q)_n=\sum_{i\ge 1}\begin{bmatrix} n\\ i\end{bmatrix}(-1)^{i-1}q^{\binom{i}{2}}x^i\,. \tag{2.13}$$

Next, we apply $\left(T_q{}^{m_k}\eta_xP_q\right)\ldots\left(T_q{}^{m_2}\eta_xP_q\right)\left(T_q{}^{m_1}\eta_xP_q\right)$ to both sides. Then, we can deduce in the same manner as above to finally arrive at (1.9) and the proof for **Theorem A** is complete.

We shall now state a corollary. Let $n_j\ge 0, 1\le j\le r$, be integers such that

$n_1+n_2+\cdots+n_r=k$. And let

$m_1=m_2=\cdots m_{n_1-1}=1, m_{n_1}=2,$

$m_{n_1+1}=m_{n_1+2}=\cdots m_{n_1+n_2-1}=1, m_{n_1+n_2}=2,$

.

.

.

$m_{n_1+\cdots+n_{r-1}+1}=m_{n_1+\cdots+n_{r-1}+2}=\cdots m_{n_1+\cdots+n_{r-1}+n_r-1}=1, m_k=1\,.$

Then, **Theorem A** becomes

**Corollary 2.2.** *For $n_j \geq 0, 1 \leq j \leq r$, there holds*

$$(a).\quad \sum_{\substack{n\geq i_1\geq\cdots\geq i_r\geq 1\\ i_0=n}} \left(q^{-i_r}-1\right)x^{i_r} \prod_{1\leq j\leq r} \frac{q^{i_j}}{\left(1-q^{i_j}\right)^2} h_{n_j-1}\left(\frac{1}{1-q^{i_j}},\dots,\frac{1}{1-q^{i_{j-1}}}\right)$$

$$=\sum_{i\geq 1}\begin{bmatrix} n\\ i\end{bmatrix}\frac{(-1)^{i-1}q^{\binom{i+1}{2}-ni}}{(1-q^i)^{n_1}} \sum_{\substack{i\geq i_1\geq\cdots\geq i_{r-1}\geq 1\\ i_0=i}} \left(1-(x;q)_{i_{r-1}}\right) \prod_{1\leq j\leq r-1} \frac{q^{i_j}}{\left(1-q^{i_j}\right)^{n_{j+1}+1}}, \tag{2.14}$$

$$(b).\quad \sum_{\substack{n\geq i_1\geq\cdots\geq i_r\geq 1\\ i_0=n}} \left(1-q^{i_r}\right)\left(1-(x;q)_{i_r}\right) \prod_{1\leq j\leq r} \frac{q^{i_j}}{\left(1-q^{i_j}\right)^2} h_{n_j-1}\left(\frac{q^{i_j}}{1-q^{i_j}},\dots,\frac{q^{i_{j-1}}}{1-q^{i_{j-1}}}\right)$$

$$=\sum_{i\geq 1}\begin{bmatrix} n\\ i\end{bmatrix}\frac{(-1)^{i-1}q^{\binom{i}{2}+in_1}}{(1-q^i)^{n_1}} \sum_{\substack{i\geq i_1\geq\cdots\geq i_{r-1}\geq 1\\ i_0=i}} x^{i_{r-1}} \prod_{1\leq j\leq r-1} \frac{q^{n_{j+1}i_j}}{\left(1-q^{i_j}\right)^{n_{j+1}+1}}. \tag{2.15}$$

It is relevant here to make some remarks as to why we use the method of Jackson integral. Why not simply use mathematical induction? Because, strictly speaking, in the context of this proof via Jackson integral, *$n$ need not to be a positive integer* in **Theorem A**. We can re-interpret this section without assuming $n$ to be a positive integer. In this sense, the proof by Jackson integral is far more powerful. Since, even by re-interpreting $(x;q)_n$ and the sum $\sum_{1\leq i\leq n} a_n$ as

$$(x;q)_n=\frac{(x;q)_\infty}{(xq^n;q)_\infty}, n\in\mathbb{C},$$

$$\sum_{1\leq i\leq n} a_n=\sum_{i\geq 1}(a_i-a_{i+n}), n\in\mathbb{C}$$

(2.5), (2.6) and (2.11) will still hold true. And consequently, **Lemma 2.1** will still hold. This is why we do not restrict $n$ in our statements. For example, re-interpreting the case $k=2, m_1=m_2=1$ of (1.9) for a non-integer $n=a$ gives us

$$\sum_{i\geq 1}\frac{(1-q^a)(1-q^{a-1})\dots\left(1-q^{a-i+1}\right)}{(1-q)(1-q^2)\dots(1-q^i)}\frac{q^{\frac{i^2+3i}{2}}\left(1-x^i\right)}{(1-q^i)^2}$$
$$=\sum_{i\geq r\geq 1}\frac{q^{i+r}}{(1-q^i)(1-q^r)}(x;q)_r-\left(\sum_{i\geq 1}\frac{q^{i+a}}{1-q^{i+a}}\right)\left(\sum_{i\geq 1}\frac{q^i}{1-q^i}(x;q)_i\right)$$
$$+\frac{(x;q)_\infty}{(xq^a;q)_\infty}\sum_{i,r\geq 1}\frac{q^{r+2i+2a}}{(1-q^{i+a})(1-q^{r+i+a})}(xq^a;q)_{r+i}. \tag{2.16}$$

Thus, it is plausible to hypothesize that **Proposition B** and **Proposition C** can also hold for a non-integer $n$. However, this question shall be left open to the readers. Our exposition shall concern with their proofs only.

## 3. Proof of Proposition B and some consequences

First, we shall prove the following lemma.

**Lemma 3.1.** *for $k \geq 0$, there holds*

$$\sum_{r\geq 0}\binom{k+r}{k}(z-1)^r h_{k+r}\left(\frac{q^i}{1-q^i},\dots,\frac{q^n}{1-q^n}\right)=\frac{(1-q^i)(q;q)_n(zq;q)_i}{(1-zq^i)(zq;q)_n(q;q)_i}h_k\left(\frac{q^i}{1-zq^i},\dots,\frac{q^n}{1-zq^n}\right). \quad (3.1)$$

*Proof of Lemma 3.1.*

Let

$$A=\sum_{k\geq 0}x^k\sum_{r\geq 0}\binom{k+r}{k}(z-1)^r h_{k+r}\left(\frac{q^i}{1-q^i},\dots,\frac{q^n}{1-q^n}\right),$$

$$B=\frac{(1-q^i)(q;q)_n(zq;q)_i}{(1-zq^i)(zq;q)_n(q;q)_i}\sum_{k\geq 0}x^k h_k\left(\frac{q^i}{1-zq^i},\dots,\frac{q^n}{1-zq^n}\right).$$

Then, it is sufficient if we prove that $A=B$.

$$\begin{aligned}A&=\sum_{k\geq 0}x^k\sum_{r\geq 0}\binom{k+r}{k}(z-1)^r h_{k+r}\left(\frac{q^i}{1-q^i},\dots,\frac{q^n}{1-q^n}\right)\\&=\sum_{m\geq 0}h_m\left(\frac{q^i}{1-q^i},\dots,\frac{q^n}{1-q^n}\right)\sum_{0\leq k\leq m}\binom{m}{k}x^k(z-1)^{m-k}=\sum_{m\geq 0}(x+z-1)^m h_m\left(\frac{q^i}{1-q^i},\dots,\frac{q^n}{1-q^n}\right)\\&=\frac{1}{\left(1-\frac{(x+z-1)q^i}{1-q^i}\right)\dots\left(1-\frac{(x+z-1)q^n}{1-q^n}\right)}=\frac{(1-q^i)\dots(1-q^n)}{(1-(x+z)q^i)\dots(1-(x+z)q^n)}.\end{aligned}$$

$$\begin{aligned}B&=\frac{(1-q^i)(q;q)_n(zq;q)_i}{(1-zq^i)(zq;q)_n(q;q)_i}\sum_{k\geq 0}x^k h_k\left(\frac{q^i}{1-zq^i},\dots,\frac{q^n}{1-zq^n}\right)\\&=\frac{(1-q^i)(q;q)_n(zq;q)_i}{(1-zq^i)(zq;q)_n(q;q)_i}\frac{1}{\left(1-\frac{xq^i}{1-zq^i}\right)\dots\left(1-\frac{xq^n}{1-zq^n}\right)}=\frac{(1-q^i)\dots(1-q^n)}{(1-(x+z)q^i)\dots(1-(x+z)q^n)}.\end{aligned}$$

Thus we have proved **Lemma 3.1**. Note that if we substitute $q^{-1}$ for $q$, we also have

$$\sum_{r\geq 0}\binom{k+r}{k}(1-z)^r h_{k+r}\left(\frac{1}{1-q^i},\dots,\frac{1}{1-q^n}\right)=z^{i-n}\frac{(1-q^i)(q;q)_n(z^{-1}q;q)_i}{(z-q^i)(z^{-1}q;q)_n(q;q)_i}h_k\left(\frac{1}{z-q^i},\dots,\frac{1}{z-q^n}\right). \quad (3.2)$$

*Proof of Proposition B.*

$.(a)$

We shall consider the expression from the left hand side. By **Lemma 3.1**, we have

$$\sum_{\substack{n\geq i_1\geq\cdots\geq i_k\geq 1\\ i_0=n}} \frac{(x;q)_{i_k}(z_kq;q)_{i_k}}{(q;q)_{i_k}} \prod_{1\leq j\leq k-1} \frac{(z_jq;q)_{i_j}}{(1-q^{i_j})(z_{j+1}q;q)_{i_j}} \prod_{1\leq j\leq k} \frac{q^{i_j}}{(1-z_jq^{i_j})} h_{l_j-1}\left(\frac{q^{i_j}}{1-z_jq^{i_j}},\dots,\frac{q^{i_{j-1}}}{1-z_jq^{i_{j-1}}}\right)$$

$$= \frac{(z_1q;q)_n}{(q;q)_n} \sum_{\substack{r_j\geq 0\\ 1\leq j\leq k}} \prod_{1\leq j\leq k} \binom{l_j-1+r_j}{l_j-1} (z_j-1)^{r_j}$$

$$\times \sum_{\substack{n\geq i_1\geq\cdots\geq i_k\geq 1\\ i_0=n}} (1-q^{i_k})(x;q)_{i_k} \prod_{1\leq j\leq k} \frac{q^{i_j}}{(1-q^{i_j})^2} h_{l_j-1+r_j}\left(\frac{q^{i_j}}{1-q^{i_j}},\dots,\frac{q^{i_{j-1}}}{1-q^{i_{j-1}}}\right).$$

And by (2.15) of **Corollary 2.2**, this equals

$$\frac{(z_1q;q)_n}{(q;q)_n} \sum_{\substack{r_j\geq 0\\ 1\leq j\leq k}} \prod_{1\leq j\leq k} \binom{l_j-1+r_j}{l_j-1} (z_j-1)^{r_j}$$

$$\times \sum_{i\geq 1} \begin{bmatrix} n\\ i\end{bmatrix} \frac{(-1)^{i-1}q^{\binom{i}{2}+i(l_1+r_1)}}{(1-q^i)^{(l_1+r_1)}} \sum_{\substack{i\geq i_1\geq\cdots\geq i_{k-1}\geq 1\\ i_0=i}} (1-x^{i_{k-1}}) \prod_{1\leq j\leq k-1} \frac{q^{(l_{j+1}+r_{j+1})i_j}}{(1-q^{i_j})^{(l_{j+1}+r_{j+1}+1)}}.$$

And since $\sum_{r\geq 0}\binom{l-1+r}{l-1}(z-1)^r \frac{q^{i(l+r)}}{(1-q^i)^{l+r}} = \frac{q^{il}}{(1-zq^i)^l}$, applying this to the above expression, we arrive at the right hand side of (1.10) and the proof is complete.

$(b)$.

In the same manner as above, with (3.2) and (2.14) in mind, we can prove (1.11).

Let $k=1, l_1=l$ and $z_1=z$, then we get the following corollary.

**Corollary 3.2.** *For $l\geq 1$, there holds*

$$(a).\ \sum_{1\leq i\leq n} \frac{(x;q)_i(zq;q)_{i-1}q^i}{(q;q)_i} h_{l-1}\left(\frac{q^i}{1-zq^i},\dots,\frac{q^n}{1-zq^n}\right) = \frac{(zq;q)_n}{(q;q)_n} \sum_{i\geq 1} \begin{bmatrix} n\\ i\end{bmatrix} \frac{(-1)^{i-1}q^{\binom{i}{2}+il}}{(1-zq^i)^l} (1-x^i), \tag{3.3}$$

$$(b).\ \sum_{1\leq i\leq n} \frac{(1-x^i)(z^{-1}q;q)_{i-1}z^{i-1}}{(q;q)_i} h_{l-1}\left(\frac{1}{z-q^i},\dots,\frac{1}{z-q^n}\right) = \frac{z^n(z^{-1}q;q)_n}{(q;q)_n} \sum_{i\geq 1} \begin{bmatrix} n\\ i\end{bmatrix} \frac{(-1)^{i-1}q^{\binom{i+1}{2}-ni}}{(z-q^i)^l} (x;q)_i. \tag{3.4}$$

We shall note again that (3.4) is identical to (1.5). If we take $l = 1$ in (3.4) and substitute $z = y^{-1}, x = w$, then we arrive at an identity due to Fu and Lascoux[12]

$$\frac{(q;q)_n}{(yq;q)_n} \sum_{1\le i\le n} \frac{w^i y^{n-i}(yq;q)_{i-1}}{(q;q)_i} = \sum_{i\ge 1} \begin{bmatrix} n \\ i \end{bmatrix} \frac{(-1)^{i-1} q^{\binom{i+1}{2}-in}}{1-yq^i} (1-(w;q)_i) \,. \tag{3.5}$$

We shall now give an identity which is an extension of (3.5). The author was unable to find this identity in q-literature.

**Proposition 3.3.** *There holds*

$$\frac{(q;q)_n}{(yz;q)_n} \sum_{1\le i\le n} \frac{w^i y^{n-i}(x;q)_i (y;q)_i (z;q)_{n-i}}{(wx;q)_i (q;q)_i (q;q)_{n-i}} = \sum_{i\ge 1} \begin{bmatrix} n \\ i \end{bmatrix} (-1)^{i-1} q^{\binom{i+1}{2}-in} \left(1-\frac{(w;q)_i}{(wx;q)_i}\right) \frac{(y;q)_i}{(yz;q)_i} \,. \tag{3.6}$$

*Proof of Proposition 3.3.*

We shall prove by means of operators. First, we apply the operator $\sum_{j\ge 0} \frac{(w;q)_j}{(q;q)_j} x^j {\eta_w}^j$ on both sides of (3.5). with $\sum_{j\ge o} \frac{(a;q)_j}{(q;q)_j} t^j = \frac{(at;q)_\infty}{(t;q)_\infty}$ (Heine[7]) in mind, we arrive at

$$\frac{(q;q)_n}{(yq;q)_n} \sum_{1\le i\le n} \frac{w^i y^{n-i}(x;q)_i (yq;q)_{i-1}}{(wx;q)_i (q;q)_i} = \sum_{i\ge 1} \begin{bmatrix} n \\ i \end{bmatrix} \frac{(-1)^{i-1} q^{\binom{i+1}{2}-in}}{1-yq^i} \left(1-\frac{(w;q)_i}{(wx;q)_i}\right) . \tag{3.7}$$

Next, apply $\sum_{j\ge 0} \frac{(-1)^j q^{\binom{j+1}{2}}}{(q;q)_j} {\eta_y}^j$ to both sides of (3.7) to get

$$\frac{(q;q)_n}{(yq;q)_n} \sum_{1\le i\le n} \frac{w^i y^{n-i}(x;q)_i (yq;q)_{i-1}}{(wx;q)_i (q;q)_i} \sum_{j\ge 0} \frac{(-1)^j q^{\binom{j+1}{2}} q^{(n-i)j} (yq^i;q)_j}{(q;q)_j (yq^{n+1};q)_j}$$
$$= \sum_{i\ge 1} \begin{bmatrix} n \\ i \end{bmatrix} (-1)^{i-1} q^{\binom{i+1}{2}-in} \left(1-\frac{(w;q)_i}{(wx;q)_i}\right) \sum_{j\ge 0} \frac{(-1)^j q^{\binom{j+1}{2}}}{(q;q)_j (1-yq^{i+j})} .$$

Since $\sum_{j\ge 0} \frac{(-1)^j q^{\binom{j+1}{2}}}{(q;q)_j (1-yq^j)} = \frac{(q;q)_\infty}{(y;q)_\infty}$ and $\sum_{j\ge 0} \frac{(-1)^j q^{\binom{j+1}{2}} q^{(n-i)j} (yq^i;q)_j}{(q;q)_j (yq^{n+1};q)_j} = \frac{(q^{n-i+1};q)_\infty}{(yq^{n+1};q)_\infty}$,(The former can be deduced from the partial fraction decomposition of $\frac{1}{(y;q)_\infty}$ , while the latter is the case $t \to 1, b = yq^n, a = q^{n-i}$ of the identity (12.3) appearing in [10,p.13]) we arrive at

$$(q;q)_n \sum_{1\le i\le n} \frac{w^i y^{n-i}(x;q)_i (y;q)_i}{(wx;q)_i (q;q)_i (q;q)_{n-i}} = \sum_{i\ge 1} \begin{bmatrix} n \\ i \end{bmatrix} (-1)^{i-1} q^{\binom{i+1}{2}-in} \left(1-\frac{(w;q)_i}{(wx;q)_i}\right) (y;q)_i \,. \tag{3.8}$$

Lastly, we use the operator $\sum_{j\ge 0} \frac{(y;q)_j}{(q;q)_j} z^j {\eta_y}^j$ on both sides of (3.8) to finally arrive at the desired result

$$\frac{(q;q)_n}{(yz;q)_n}\sum_{1\le i\le n}\frac{w^i y^{n-i}(x;q)_i(y;q)_i(z;q)_{n-i}}{(wx;q)_i(q;q)_i(q;q)_{n-i}}=\sum_{i\ge 1}\begin{bmatrix} n\\ i\end{bmatrix}(-1)^{i-1}q^{\binom{i+1}{2}-in}\left(1-\frac{(w;q)_i}{(wx;q)_i}\right)\frac{(y;q)_i}{(yz;q)_i}.$$

Let $x=0$, then we have

$$\frac{(q;q)_n}{(yz;q)_n}\sum_{1\le i\le n}\frac{w^i y^{n-i}(y;q)_i(z;q)_{n-i}}{(q;q)_i(q;q)_{n-i}}=\sum_{i\ge 1}\begin{bmatrix} n\\ i\end{bmatrix}(-1)^{i-1}q^{\binom{i+1}{2}-in}(1-(w;q)_i)\frac{(y;q)_i}{(yz;q)_i}. \tag{3.9}$$

Equating the coefficients of $w^r$ from both sides, we get the following corollary.

**Corollary 3.4.** *There holds*

$$\sum_{i\ge 1}\begin{bmatrix} n\\ i\end{bmatrix}\begin{bmatrix} i\\ r\end{bmatrix}(-1)^{i-1}q^{\binom{i+1}{2}-in}\frac{(y;q)_i}{(yz;q)_i}=y^{n-r}(-1)^{r-1}q^{-\binom{r}{2}}\frac{(q;q)_n(y;q)_r(z;q)_{n-r}}{(yz;q)_n(q;q)_r(q;q)_{n-r}}. \tag{3.10}$$

## 4. Proof of Proposition C

Our proof of **Proposition C** requires a preliminary result to be established beforehand.

**Proposition 4.1.** *For k arbitrary constants $z_1,\dots,z_k$, there holds*

$$\sum_{i_1=1}^{n}\frac{q^{i_1}}{(1-q^{i_1})(1-z_1q^{i_1})}\sum_{i_2=1}^{i_1}\frac{q^{i_2}}{(1-q^{i_2})(1-z_2q^{i_2})}\cdots\sum_{i_k=1}^{i_{k-1}}\frac{q^{i_k}}{(1-q^{i_k})(1-z_kq^{i_k})}(1-x^{i_k})$$

$$=\sum_{i\ge 1}\begin{bmatrix} n\\ i\end{bmatrix}(-1)^{i-1}q^{\binom{i+1}{2}-ni}\frac{(q;q)_{i-1}}{(z_1q;q)_i}$$

$$\times\sum_{i\ge i_1\ge\cdots\ge i_k\ge 1}\frac{(x;q)_{i_k}(z_kq;q)_{i_k}}{(q;q)_{i_k}}\prod_{1\le j\le k-1}\frac{(z_jq;q)_{i_j}}{(1-q^{i_j})(z_{j+1}q;q)_{i_j}}\prod_{1\le j\le k}\frac{q^{i_j}}{1-z_jq^{i_j}}. \tag{4.1}$$

*Proof of Proposition 4.1.*

We shall proceed from the left hand side.

$$\sum_{i_1=1}^{n}\frac{q^{i_1}}{(1-q^{i_1})(1-z_1q^{i_1})}\sum_{i_2=1}^{i_1}\frac{q^{i_2}}{(1-q^{i_2})(1-z_2q^{i_2})}\cdots\sum_{i_k=1}^{i_{k-1}}\frac{q^{i_k}}{(1-q^{i_k})(1-z_kq^{i_k})}(1-x^{i_k})$$

$$=\sum_{\substack{m_j\ge 1\\ 1\le j\le k}}\prod_{1\le j\le k}(z_j-1)^{m_j-1}\sum_{n\ge i_1\ge\cdots\ge i_k\ge 1}(1-x^{i_k})\prod_{1\le j\le k}\frac{q^{m_ji_j}}{(1-q^{i_j})^{m_j+1}}$$

$$= \sum_{\substack{m_j \geq 1 \\ 1 \leq j \leq k}} \prod_{1 \leq j \leq k} (z_j - 1)^{m_j - 1} \sum_{i \geq 1} \begin{bmatrix} n \\ i \end{bmatrix} \frac{(-1)^{i-1} q^{\binom{i+1}{2} - ni}}{1 - q^i}$$

$$\times \sum_{\substack{i \geq i_1 \geq \cdots \geq i_k \geq 1 \\ i_0 = i}} (1 - q^{i_k})(x;q)_{i_k} \prod_{1 \leq j \leq k} \frac{q^{i_j}}{(1 - q^{i_j})^2} h_{m_j - 1}\left(\frac{q^{i_j}}{1 - q^{i_j}}, \ldots, \frac{q^{i_{j-1}}}{1 - q^{i_{j-1}}}\right)$$

$$= \sum_{i \geq 1} \begin{bmatrix} n \\ i \end{bmatrix} (-1)^{i-1} q^{\binom{i+1}{2} - ni} \frac{(q;q)_{i-1}}{(z_1 q;q)_i}$$

$$\times \sum_{i \geq i_1 \geq \cdots \geq i_k \geq 1} \frac{(x;q)_{i_k} (z_k q;q)_{i_k}}{(q;q)_{i_k}} \prod_{1 \leq j \leq k-1} \frac{(z_j q;q)_{i_j}}{(1 - q^{i_j})(z_{j+1} q;q)_{i_j}} \prod_{1 \leq j \leq k} \frac{q^{i_j}}{1 - z_j q^{i_j}}.$$

Thus, our proof is complete. Now we shall proceed to prove **Proposition C**.

*Proof of Proposition C.*

Let $z_j = zq^{-j}$ in **Proposition 4.1** to get

$$\sum_{i_1=1}^{n} \frac{q^{i_1}}{(1 - q^{i_1})(1 - zq^{i_1 - 1})} \sum_{i_2=1}^{i_1} \frac{q^{i_2}}{(1 - q^{i_2})(1 - zq^{i_2 - 2})} \cdots \sum_{i_k=1}^{i_{k-1}} \frac{q^{i_k}}{(1 - q^{i_k})(1 - zq^{i_k - k})} (1 - x^{i_k})$$

$$= \frac{1}{(zq^{-k};q)_k} \sum_{i \geq 1} \begin{bmatrix} n \\ i \end{bmatrix} (-1)^{i-1} q^{\binom{i+1}{2} - ni} \frac{(q;q)_{i-1}}{(z;q)_i} \sum_{i \geq i_1 \geq \cdots \geq i_k \geq 1} \frac{(x;q)_{i_k} (zq^{-k};q)_{i_k}}{(q;q)_{i_k - 1}} \prod_{1 \leq j \leq k} \frac{q^{i_j}}{(1 - q^{i_j})}. \qquad (4.2)$$

Then apply the operator $\sum_{j \geq 0} \frac{(x;q)_j}{(q;q)_j} t^j \eta_x{}^j$ on both sides to arrive at

$$\sum_{i_1=1}^{n} \frac{q^{i_1}}{(1 - q^{i_1})(1 - zq^{i_1 - 1})} \sum_{i_2=1}^{i_1} \frac{q^{i_2}}{(1 - q^{i_2})(1 - zq^{i_2 - 2})} \cdots \sum_{i_k=1}^{i_{k-1}} \frac{q^{i_k}}{(1 - q^{i_k})(1 - zq^{i_k - k})} \left(1 - x^{i_k} \frac{(t;q)_{i_k}}{(xt;q)_{i_k}}\right)$$

$$= \frac{1}{(zq^{-k};q)_k} \sum_{i \geq 1} \begin{bmatrix} n \\ i \end{bmatrix} (-1)^{i-1} q^{\binom{i+1}{2} - ni} \frac{(q;q)_{i-1}}{(z;q)_i} \sum_{i \geq i_1 \geq \cdots \geq i_k \geq 1} \frac{(x;q)_{i_k} (zq^{-k};q)_{i_k}}{(xt;q)_{i_k} (q;q)_{i_k - 1}} \prod_{1 \leq j \leq k} \frac{q^{i_j}}{(1 - q^{i_j})}. \qquad (4.3)$$

Divide both sides by $1 - x$ and let $x \to 1$. Since $\lim_{x \to 1} \frac{1 - x^{i_k} \frac{(t;q)_{i_k}}{(xt;q)_{i_k}}}{1 - x} = \sum_{1 \leq r \leq i_k} \frac{1}{1 - tq^{r-1}}$, we have

$$\sum_{i_1=1}^{n} \frac{q^{i_1}}{(1 - q^{i_1})(1 - zq^{i_1 - 1})} \sum_{i_2=1}^{i_1} \frac{q^{i_2}}{(1 - q^{i_2})(1 - zq^{i_2 - 2})} \cdots \sum_{i_k=1}^{i_{k-1}} \frac{q^{i_k}}{(1 - q^{i_k})(1 - zq^{i_k - k})} \sum_{r=1}^{i_k} \frac{1}{1 - tq^{r-1}}$$

$$= \frac{1}{(zq^{-k};q)_k} \sum_{i \geq 1} \begin{bmatrix} n \\ i \end{bmatrix} (-1)^{i-1} q^{\binom{i+1}{2} - ni} \frac{(q;q)_{i-1}}{(z;q)_i} \sum_{i \geq i_1 \geq \cdots \geq i_k \geq 1} \frac{(zq^{-k};q)_{i_k}}{(t;q)_{i_k}} \prod_{1 \leq j \leq k} \frac{q^{i_j}}{(1 - q^{i_j})}. \qquad (4.4)$$

And since

$$\sum_{\substack{i\geq i_1\geq\cdots\geq i_k\geq 1\\ i_0=i}} \frac{(zq^{-k};q)_{i_k}}{(t;q)_{i_k}} \prod_{1\leq j\leq k} \frac{q^{i_j}}{(1-q^{i_j})}$$

$$= \frac{(tz^{-1}q^k;q)_\infty}{(t;q)_\infty} \left( \sum_{j\geq 0} \frac{(zq^{-k};q)_j}{(q;q)_j} (tz^{-1}q^k)^j {\eta_z}^j \right) \sum_{i\geq i_1\geq\cdots\geq i_k\geq 1} (zq^{-k};q)_{i_k} \prod_{1\leq j\leq k} \frac{q^{i_j}}{1-q^{i_j}}$$

$$= \frac{(tz^{-1}q^k;q)_\infty}{(t;q)_\infty} \left( \sum_{j\geq 0} \frac{(zq^{-k};q)_j}{(q;q)_j} (tz^{-1}q^k)^j {\eta_z}^j \right) \sum_{r\geq 1} \begin{bmatrix} i \\ r \end{bmatrix} \frac{(-1)^{r-1} q^{\binom{r}{2}+rk}}{(1-q^r)^k} (1 - z^r q^{-rk})$$

$$= \sum_{r\geq 1} \begin{bmatrix} i \\ r \end{bmatrix} \frac{(-1)^{r-1} q^{\binom{r}{2}+rk}}{(1-q^r)^k} \left( 1 - z^r q^{-rk} \frac{(tz^{-1}q^k;q)_r}{(t;q)_r} \right), \tag{4.5}$$

(4.4) is

$$\sum_{i_1=1}^{n} \frac{q^{i_1}}{(1-q^{i_1})(1-zq^{i_1-1})} \sum_{i_2=1}^{i_1} \frac{q^{i_2}}{(1-q^{i_2})(1-zq^{i_2-2})} \cdots \sum_{i_k=1}^{i_{k-1}} \frac{q^{i_k}}{(1-q^{i_k})(1-zq^{i_k-k})} \sum_{r=1}^{i_k} \frac{1}{1-tq^{r-1}}$$

$$= \frac{1}{(zq^{-k};q)_k} \sum_{r\geq 1} \frac{(-1)^{r-1} q^{\binom{r}{2}+rk}}{(1-q^r)^k} \left( 1 - z^r q^{-rk} \frac{(tz^{-1}q^k;q)_r}{(t;q)_r} \right) \sum_{i\geq 1} \begin{bmatrix} n \\ i \end{bmatrix} \begin{bmatrix} i \\ r \end{bmatrix} (-1)^{i-1} q^{\binom{i+1}{2}-ni} \frac{(q;q)_{i-1}}{(z;q)_i}. \tag{4.6}$$

To evaluate the inner sum from the right hand side, we divide both sides of **Corollary 3.4** by $1-y$ and set $y=1$ to get

$$\sum_{i\geq 1} \begin{bmatrix} n \\ i \end{bmatrix} \begin{bmatrix} i \\ r \end{bmatrix} (-1)^{i-1} q^{\binom{i+1}{2}-in} \frac{(q;q)_{i-1}}{(z;q)_i} = \frac{(-1)^{r-1} q^{-\binom{r}{2}}}{1-q^r} \frac{(q;q)_n (z;q)_{n-r}}{(z;q)_n (q;q)_{n-r}}. \tag{4.7}$$

Thus, (4.6) can finally be transformed into

$$\sum_{i_1=1}^{n} \frac{q^{i_1}}{(1-q^{i_1})(1-zq^{i_1-1})} \sum_{i_2=1}^{i_1} \frac{q^{i_2}}{(1-q^{i_2})(1-zq^{i_2-2})} \cdots \sum_{i_k=1}^{i_{k-1}} \frac{q^{i_k}}{(1-q^{i_k})(1-zq^{i_k-k})} \sum_{r=1}^{i_k} \frac{1}{1-tq^{r-1}}$$

$$= \sum_{r\geq 1} \begin{bmatrix} n \\ r \end{bmatrix} \frac{(q;q)_{r-1} (z;q)_{n-r}}{(zq^{-k};q)_{n+k} (1-q^r)^k} \left( 1 - z^r q^{-rk} \frac{(tz^{-1}q^k;q)_r}{(t;q)_r} \right) q^{rk},$$

and the proof is complete.

We now end our exposition with the following corollary. Let $n \to \infty$ in **Proposition C**, to get

**Corollary 4.2.** *For natural number k, there holds*

$$\sum_{i_1=1}^{\infty}\frac{q^{i_1}}{(1-q^{i_1})(1-zq^{i_1-1})}\sum_{i_2=1}^{i_1}\frac{q^{i_2}}{(1-q^{i_2})(1-zq^{i_2-2})}\cdots\sum_{i_k=1}^{i_{k-1}}\frac{q^{i_k}}{(1-q^{i_k})(1-zq^{i_k-k})}\sum_{r=1}^{i_k}\left(\frac{1}{1-zq^{r-k-1}}-\frac{1}{1-tq^{r-1}}\right)$$

$$=\frac{1}{(zq^{-k};q)_k}\sum_{r=1}^{\infty}\frac{(tz^{-1}q^k;q)_r z^r}{(t;q)_r(1-q^r)^{k+1}}\,. \qquad (4.8)$$

## *References*

[1] H. Prodinger, "A q-ANALOGUE OF A FORMULA OF HERNANDEZ OBTAINED BY INVERTING A RESULT OF DILCHER," 1999.

[2] H. Prodinger, "q-IDENTITIES OF FU AND LASCOUX PROVED BY THE q-RICE FORMULA," 2004.

[3] V. Hernández, "Solution IV(A Reciprocal Summation Identity: 10490)," *The American Mathematical Monthly,* vol. 106, p. 589, 1999.

[4] K. Dilcher, "Some q-series identities related to divisor functions," *Discrete Mathematics,* vol. 145, pp. 83-93, 1995.

[5] J. Kluyver, "Vraagstuk XXXVII. (Solution by S.C. van Veen)," *Wiskundige Opgaven,* pp. 92-93, 1919.

[6] F. H. Jackson, "On q-Definite Integrals," *The Quarterly Journal of Pure and Applied Mathematics,* vol. 41, pp. 193-203, 1910.

[7] E. Heine, "Untersuchungen über die Reihe...," *J. Reine Angew. Math,* vol. 34, pp. 285-328, 1847.

[8] G. E. Andrews, "On the Foundations of Combinatorial Theory V, Eulerian Differential Operators," *Studies in Applied Mathematics,* vol. L, pp. 345-375, 1971.

[9] J. Zeng, "On some q-identities related to divisor functions," *Advances in Applied Mathematics,* vol. 34, pp. 313-315, 2005.

[10] N. J. Fine, BASIC HYPERGEOMETRIC SERIES AND APPLICATIONS, American Mathematical Society, 1988.

[11] V. J. W. Guo and C. Zhang, "Some further q-series identities related to divisor functions," *The Ramanujan Journal,* 2011.

[12] A. M. Fu and A. Lascoux, "q-Identities related to overpartitions and divisor functions," *THE ELECTRONIC JOURNAL OF COMBINATORICS,* 2005.

[13] A. Xu and Z. Cen, "Combinatorial identities from contour integrals of rational functions," *The Ramanujan Journal,* pp. 103-114, 2015.

[14] A. Xu, "On a general q-identity," *THE ELECTRONIC JOURNAL OF COMBINATORICS,* 2014.

[15] D. M. Bradley, "Duality for finite multiple harmonic q-series," *DISCRETE MATHEMATICS,* pp. 44-56, 2005.

[16] V. J. Guo and J. Zeng, "Basic and bibasic identities related to divisor functions," *Journal of Mathematical Analysis and Applications,* pp. 1197-1209, 2015.

[17] M. E. Ishmail and D. Stanton, "Some Combinatorial and Analytical Identities," *Annals of Combinatorics,* 2010.

[18] V. Kac and P. Cheung, Quantum Calculus, Springer-Verlag New York, Inc., 2002.